 \documentclass[12pt,leqno]{amsart}
\newtheorem{dummy}{}[section]
\newtheorem{definition}[dummy]{Definition}
\newtheorem{theorem}[dummy]{Theorem}

\newtheorem{proposition}[dummy]{Proposition}
\newtheorem{lemma}[dummy]{Lemma}

\newtheorem{example}[dummy]{Example}
\newtheorem{question}[dummy]{Question}
\begin{document}
\bibliographystyle{plain}
\title{Derived Langlands VII: The PSH algebra of products of general linear groups } 
\author{Victor P. Snaith}
\date{16 December 2020}

\maketitle
 \tableofcontents 
 
{\it My father, Victor Snaith, died on July 3, 2021 after a battle with a blood disorder named myelofibrosis. In the final weeks and months of his life he completed Derived Langlands VIII published it online. Previous to that he had written this, the seventh part of his Derived Langlands series but, as time was running short and his health worsenend, had not published it, despite having mentioned to me in conversation that it was essentially complete and that I was welcome to publish it for him posthumously. I do not know whether he had not published it because there were still details to attend to, so this document is best treated as unfinished. I have made notes in the text below, labelled DS, where he had left remarks in the typesetting that he presumably would have intended to remove before publication if he had more time.}

{\it Regards,
Dan Snaith}

 \section{Local fields remarks}
 
(i)  \   I am going to spend a lot of space and effort putting a very elaborate PSH-like  structure on the $R_{+}(-)$ groups of products of finite general linear groups. This is not the case we want. Firstly one would really want the actual big PSH algebra of products of general linear groups with entries in a characteristic zero  $p$-adic local field. There may be technical difficulties wth this. However the $R_{+}(-)$ gadget for products of general linear groups with entries in a characteristic zero  $p$-adic local field seems to work for us by allowing various reduction to compact open subgroups and reduction maps modulo different prime powers from there. These reductions may allow the verification of functional equations and analytic groups properties which characterise the Riemann zeta function and presumably similarly characterise the $2$-variable L-functions.

{ii) \  The monomial resolution in the case of admissible representations (in the various senses of 
(\cite{Sn18}, \cite{Sn20}, \cite{Sn20b}, \cite{newindnotes}, \cite{Sn20c}, \cite{Sn20d})
of general linear groups of local fields
(or more generally locally $p$-adic reductive algebraic groups) is, inter alia, a Grothendieck-style approach to Langlands Quotient Theorem. This is one of my favourite classification theorems
(\cite{RPL}, \cite{AJS78}, \cite{AJS79},   \cite{LLProb1970}) and in the original treatment it was a very delicate application of the methods of Harish-Chandra \cite{HC66} which proves that every admissible irreducible $V$ is a quotient of an admissible representation induced from one of a parabolic subgroup, classifying up to an equivalence relation the data which gives rise to the specific admissible irreducible under consideration.

The  Grothendieck-style approach is to manufacture a chain homotopy category, with as simple as possible adjectival defining properties (in our case the monomial category), with a resolution of the representation by a 
chain complex which is ``exact'' in a suitably enriched sense. The classifying surjection is then a (monomial category) homotopy class of the surjection from a monomial object onto $V$.

\begin{question}{$_{}$}
\label{A1.1}
\begin{em}

The monomial objects in the canonical homotopy class of monomial surjections is, as a representation, a sum of representations induced from characters of compact open modulo the centre subgroups of, for example, products of general linear groups of a local field. In order to benefit from the inductive property of Langlands theorem, which reduces to parabolics associated to products of smaller general linear groups, I would like to establish a result which showed that the monomial surjection was homotopy equivalent to one associated to products of smaller general linear groups.

Is such a result even possible?
\end{em}
\end{question}

{\bf Note} \  This is going to be another rather strange\footnote{The strangeness is a consequence of a medical time-table beyond my control \cite{JMOS18}. I shall try my best to assemble enough of a coherent further bunch of slightly non-classical Langlands notions, enough perhaps to interest the occasional reader.} and sketchy 
essay in my Derived Langlands series previously consisting of a research monograph and four essays (\cite{Sn18}, \cite{Sn20}, \cite{Sn20b}, \cite{newindnotes}, \cite{Sn20c}, \cite{Sn20d}). {\it DS: shortly before his death my dad added Derived Langlands VIII to the arXiv \cite{Sn21}.}

 \section{Positive Selfadjoint Hopf algebras}
 
Connected PSH algebras over ${\mathbb Z}$ are defined and analysed in detail in \cite{AVZ81} (see also \cite{Sn18} Chapter Nine and Appendix Three). If $H$ is such an algebra we have a multiplication
\[ m: H \otimes_{{\mathbb Z}} H = H \otimes H   \longrightarrow H \]
and a comultiplication
\[ m^{*} : H \longrightarrow H \otimes H .  \]
If $(H_{1}, m_{1}, m_{1}^{*} )$ and $(H_{2}, m_{2}, m_{2}^{*})$ are two PSH algebras then $H_{1} \otimes H_{2}$ is also a PSH algebra with multiplication
\[ \hat{m} :  (H_{1} \otimes H_{2})  \otimes (H_{1} \otimes H_{2}) \longrightarrow  H_{1} \otimes H_{2}  \]
given by the composition
\[ H_{1} \otimes H_{2}  \otimes H_{1} \otimes H_{2}  \stackrel{1 \otimes T \otimes 1}{\longrightarrow}  
H_{1} \otimes H_{1}  \otimes H_{2} \otimes H_{2}  \stackrel{m_{1} \otimes m_{2}}{\longrightarrow}  
 H_{1} \otimes H_{2}  \]  
and comultiplication
\[  \hat{m}^{*} : (H_{1} \otimes H_{2}) \longrightarrow  (H_{1} \otimes H_{2})  \otimes (H_{1} \otimes H_{2}) \]
given by the composition
\[ H_{1} \otimes H_{2} \stackrel{ m_{1}^{*} \otimes m_{2}^{*} }{\longrightarrow}  H_{1} \otimes H_{1}  \otimes H_{2} \otimes H_{2}  \stackrel{ 1 \otimes T \otimes 1}{\longrightarrow}  H_{1} \otimes H_{2}  \otimes H_{1} \otimes H_{2} . \]
The Hopf algebra conditions satisfied by $m$ and $m^{*}$ is that $m$ is a map of coalgebras and $m^{*}$ is a map of algebras.

Choose a finite fields ${\mathbb F}_{q}$ and write $G_{a} = GL_{a}{\mathbb F}_{q}$ for 
$a = 0, 1, 2, 3, \ldots$\footnote{Sometimes alternatively written $G(a)$}.
 Write $R_{a} $ for the complex representation ring $R(G_{a})$ with the interpretation that $R(G_{0}) = {\mathbb Z}$.
 
 The graded ring
 \[ R_{1} = R(0) \oplus R(1) \oplus  R(2) \oplus \ldots \]
 is a ${\mathbb Z}$-graded PSH algebra \cite{AVZ81} (see also \cite{Sn18} Chapter Nine). The graded ring
 \[ \begin{array}{ll}
 R_{2} &= (R(0) \otimes R(0)) \\
 \\
 &\oplus (R(1) \otimes R(0)) \oplus (R(0) \otimes R(1)) \\
 \\
 & \oplus
 (R(2) \otimes R(0)) \oplus (R(1) \otimes R(1)) \oplus (R(0) \otimes R(2)) \\
 \\
 & \oplus  (R(3) \otimes R(0)) \oplus (R(2) \otimes R(1)) \oplus (R(1) \otimes R(2)) \oplus (R(0) \otimes R(3) \\
 \\
 & \ldots \hspace{20pt} \ldots   \hspace{20pt} \ldots  
 \end{array}  \]
 is the tensor product of the PSH algebra $R_{1}$ with itself. It is graded by the set of ordered pairs of non-negative integers.
 
 Similarly we have the $m$-fold tensor product of the PSH algebra $R_{1}$
 \[ R_{n} = \ldots \oplus (R(a_{1}) \otimes R(a_{2}) \otimes R(a_{3})  \otimes \ldots \otimes  R(a_{n-1}) \otimes R(a_{n})) \oplus \ldots  \]
 which is graded by ordered $n$-tuples of non-negative integers.
 
 The multiplication on $R_{n}$ is given by 
  \[ m: (R(a_{1}) \otimes \dots \ \otimes R(a_{n})) \otimes  (R(b_{1}) \otimes \dots \ \otimes R(b_{n}))
  \longrightarrow (R(a_{1}+ b_{1}) \otimes \dots \ \otimes R(a_{n}+b_{n}))\]
 sending
 $(V_{1} \otimes \ldots \otimes V_{n}) \otimes (W_{1} \otimes \ldots \otimes W_{n}) $ to
 \[   {\rm Ind}_{P_{a_{1} , b_{1}}}^{G_{a_{1}+b_{1}}}( {\rm Inf}_{G_{a_{1}} \times G_{b_{1}} }^{P_{a_{1} , b_{1}}}(V_{1} \otimes W_{1}))  \otimes \ldots \otimes  {\rm Ind}_{P_{a_{n} , b_{n}}}^{G_{a_{n}+b_{n}}}( {\rm Inf}_{G_{a_{n}} \times G_{b_{n}} }^{P_{a_{n} , b_{n}}}(V_{n} \otimes W_{n}))  . \]
 
 The $(\alpha_{1} , \dots , \alpha_{n} ,  \beta_{1} , \dots , \beta_{n})$-component of the comultiplication
 \linebreak
 $m^{*}(V_{1} \otimes \ldots \otimes V_{n})$ is given by
 \[ \tau(   m^{*}(V_{1}) \otimes   m^{*}(V_{2} \otimes  \ldots \otimes m^{*}(V_{n}    ))   \]
 where $\tau$ is the shuffle given by the fermutation
 \[ (1,2,\ldots,n,n+1, n+2, \ldots , 2n) \stackrel{\tau}{\mapsto}  (1,n+1,2, n+2, \ldots,n-1,2n-1, n, 2n)\]
 
Each of the maps which inserts $k$ zeros into $(a_{1} \otimes \dots \ \otimes a_{n})$ induces a PSH inclusion of $R_{n}$ into $R_{n+k}$ 
 so that we have a big PSH algebra given by $R_{\infty} = \bigcup_{n \geq 1} R_{n}$.
 
 Since $R(G \times J) \cong R(G) \otimes R(J)$ we may rewrite $R_{n}$  for $0  \leq n \leq \infty$ 
 as 
 \[  R_{n} = \oplus_{a_{1}, \ldots , a_{n}}  \ R(G(a_{1}) \times \ldots \times G(a_{n})) .\]
 
 Now consider
  \[  R_{+, n} = \oplus_{a_{1}, \ldots , a_{n}}  \ R_{+}(G(a_{1}) \times \ldots \times G(a_{n})) \]
  where $R_{+}(G)$ is defined in (\cite{Sn94} Definition 2.2.1 p. 32) to be the free abelian group on $G$-conjugacy classes of characters $\phi : H \longrightarrow {\mathbb C}^{*}$ were $H \leq G$.
  
  Each of the pairs $(H, \phi)$ above will satisfy conditions involving a ``central character'' $\underline{\phi} : Z(G) \longrightarrow  k^{*}$ where, in the case of $G$ being a product of $GL_{m}{\mathbb F}_{q}$'s the centre will be ${\mathbb F}_{q}^{*}$ in each case. We shall assume that $Z(G) \leq H$ and that $\phi$ coincides with  $\underline{\phi}$ on $Z(G)$.

I am going to go into lots of PSH and PSH-like details in the case of products of general linear groups of 
 ${\mathbb F}_{q}$ for which example we could omit reference to  a ``central character''. However, the PSH-like structure which I shall describe on $R_{+}(-)$ in the finite field case seems, I believe, to generalise to the case of $p$-adic local fields, in which case - as in \cite{Sn18} - the treatment parcels admissible representations together according to their central character.
  
   We denote the class of $(H, \phi)$ by $(H, \phi)^{G} \in R_{+}(G)$. $R_{+}(G)$ is a ring whose multiplication is defined by a Double Coset Formula which is given in (\cite{Sn94} Exercise 2.5.7 p.68) and is given explicitly by the formula
  \[  (K, \phi)^{G} \cdot  (H, \psi)^{G} = \sum_{w \in K \backslash G / H} \   ((w^{-1}Kw) \bigcap H, w^{*}(\phi) \psi )^{G} ,\]
  where $w^{*}(\phi) (z) = \phi(wzw^{-1})$ for $z \in w^{-1}Kw$.
  
  There is a ring homomorphism
  \[   b : R_{+}(G) \longrightarrow R(G) \]
  given by $b((H, \phi)^{G}) = {\rm Ind}_{H}^{G}(\phi)$.
  
  Now let us examine $b_{n} = b : R_{+,n} \longrightarrow R_{n}$ for $0  \leq n \leq \infty$, which is an algebra map.

  \section{  The case when the central character is trivial.}

  Consider $(H, \phi)^{G(a)} \in R_{+}(G(a)) \subset R_{+1,1}$. The $(\alpha, a - \alpha)$-component of
  $m^{*}(b((H, \phi)^{G(a)})))$ is the $G(\alpha) \times G(a - \alpha)$-representation
\[ {\rm Ind}_{H}^{G(a) }( \phi))^{U_{\alpha, a - \alpha} } =   \sum_{t \in U_{\alpha, a - \alpha}  \backslash G(a) / H } \ ({\rm Ind}_{U_{\alpha, a - \alpha}  \cap tHt^{-1} }^{U_{\alpha, a - \alpha} }((t^{-1})^{*}(\phi)))^{U_{\alpha, a - \alpha} } .\]
The invariants $({\rm Ind}_{U_{\alpha, a - \alpha}  \cap tHt^{-1} }^{U_{\alpha, a - \alpha} }((t^{-1})^{*}(\phi)))^{U_{\alpha, a - \alpha} }$ are one-dimensional with a basis given by 
\[ \frac{1}{| U_{\alpha, a - \alpha}  |} \Sigma_{u \in U_{\alpha, a - \alpha}} \ u \otimes_{U_{\alpha, a - \alpha}  \cap tHt^{-1}} 1 .\]

The Double coset Formula (\cite{Sn18} p.186; \cite{Sn94} Theorem 1.2.40) is an isomorphism
\[ \alpha : {\rm Res}_{J}^{G} {\rm Ind}_{H}^{G}(\phi) \stackrel{\cong}{\longrightarrow} 
\Sigma_{z \in J \backslash G / H} \ {\rm Ind}_{J \cap zHz^{-1}}^{J}((z^{-1})^{*}(\phi)) \]
given by $\alpha( g \otimes_{H} w) = j \otimes_{J \cap zHz^{-1}}  hw$ for $g = jzh, j \in J, h \in H$ so that $\alpha^{-1}( j \otimes_{J \cap zHz^{-1}} w) = jz \otimes_{H} w$. Therefore, if $u \in U_{\alpha, a - \alpha}, h \in H$ and $ g = uth$, then $u \otimes_{U_{\alpha, a - \alpha}  \cap tHt^{-1} } 1$ corresponds to $g \otimes_{H} \phi(h)^{-1} = ut \otimes_{H} 1 \in  {\rm Ind}_{H}^{G(a)}(\phi) $.

We can simplify the expression 
\[  \begin{array}{l}
\frac{1}{| U_{\alpha, a - \alpha}  |} \Sigma_{u \in U_{\alpha, a - \alpha}} \ u \otimes_{U_{\alpha, a - \alpha}  \cap tHt^{-1}} 1 \\
\\
= \frac{1}{| U_{\alpha, a - \alpha}  |} \sum_{s=1}^{[U_{\alpha, a - \alpha} : U_{\alpha, a - \alpha} \cap tHt^{-1}]}
u_{s} \otimes_{U_{\alpha, a - \alpha}  \cap tHt^{-1} } ( \sum_{u' \in  U_{\alpha, a - \alpha} \cap tHt^{-1}} \phi(t^{-1} u' t) ) ,
\end{array} \]
where $u_{1}, u_{2}, \ldots $ are coset representatives of $U_{\alpha, a - \alpha}/ U_{\alpha, a - \alpha} \bigcap tHt^{-1}$.
However $\sum_{u' \in  U_{\alpha, a - \alpha} \cap tHt^{-1}} \phi(t^{-1} u' t) =0$ unless $\phi$ is trivial on $t^{-1}U_{\alpha, a - \alpha}t \bigcap H$ in which case 
$\sum_{u' \in  U_{\alpha, a - \alpha} \cap tHt^{-1}} \phi(t^{-1} u' t)
 = | t^{-1}U_{\alpha, a - \alpha}t \bigcap H |$\footnote{Note that for a double coset $t \in U_{\alpha, a - \alpha} \backslash G /H$ the condition that $\phi = 1$ on $t^{-1}  U_{\alpha, a - \alpha}  t \bigcap H$ holds if and only if $\phi = 1$ on  $(u'th)^{-1}  U_{\alpha, a - \alpha}  (u'th) \bigcap H = h^{-1} (t^{-1}  U_{\alpha, a - \alpha}  t \bigcap H )h$ for $u' \in U, h \in H$. }

Therefore we obtain
\[  \begin{array}{l}
\frac{1}{| U_{\alpha, a - \alpha}  |} \Sigma_{u \in U_{\alpha, a - \alpha}} \ u \otimes_{U_{\alpha, a - \alpha}  \cap tHt^{-1}} 1 \\
\\
= \left\{  \begin{array}{l}
\frac{1}{[U_{\alpha, a - \alpha} : U_{\alpha, a - \alpha} \cap tHt^{-1}]}
\sum_{s=1}^{[U_{\alpha, a - \alpha} : U_{\alpha, a - \alpha} \cap tHt^{-1}]}  u_{s} \otimes_{U_{\alpha, a - \alpha}  \cap tHt^{-1} } 1 \\
\\
\hspace{50pt}  {\rm if \ \phi = 1 \   on \ tU_{\alpha, a - \alpha}t^{-1} \bigcap H}  \\
\\
0 \  {\rm otherwise} .
\end{array}  \right.
 \end{array} \]
 From the formula for the map $\alpha^{-1}$ if $g = u_{s} t h_{s}$ with $h_{s} \in H$ then $u \otimes_{U_{\alpha, a - \alpha}  \cap tHt^{-1} } 1$ corresponds to $g \otimes_{H} \phi(h_{s})^{-1} = u_{s}t \otimes_{H} 1 \in  {\rm Ind}_{H}^{G(a)}(\phi) $ and
 \[  \begin{array}{l}
\frac{1}{| U_{\alpha, a - \alpha}  |} \Sigma_{u \in U_{\alpha, a - \alpha}} \ u \otimes_{U_{\alpha, a - \alpha}  \cap tHt^{-1}} 1 \\
\\
= \left\{  \begin{array}{l}
\frac{1}{[U_{\alpha, a - \alpha} : U_{\alpha, a - \alpha} \cap tHt^{-1}]}
\sum_{s=1}^{[U_{\alpha, a - \alpha} : U_{\alpha, a - \alpha} \cap tHt^{-1}]}  u_{s}t \otimes_{H} 1 \\
\\
\hspace{50pt}  {\rm if \  g = u_{s} t h_{s} \ and \  \phi = 1 \   on \ tU_{\alpha, a - \alpha}t^{-1} \bigcap H}  \\
\\
0 \  {\rm otherwise} .
\end{array}  \right.
 \end{array} \]
 
 Now consider the action of $G_{\alpha} \times G(a - \alpha) \cong  P_{\alpha, , a - \alpha}/U_{\alpha, , a - \alpha}$ on this element.
 
 If $r \in P_{\alpha, a - \alpha}$ and $\hat{u} \in U_{\alpha, a - \alpha}$ the action is given by
 \[  ( r \hat{u} \cdot -): g \otimes_{H} \phi(h_{s})^{-1} = u_{s}t \otimes_{H} 1  \mapsto r \hat{u} g \otimes_{H} \phi(h_{s})^{-1} =  r \hat{u} u_{s}t \otimes_{H} 1   \]
 and $ r \hat{u}g =  r \hat{u}  u_{s} t h_{s}$. Also $\hat{u} u_{s} = u_{\sigma(s)} \tilde{u}$ for some 
 $ \tilde{u} \in  U_{\alpha, a - \alpha} \bigcap tHt^{-1}$. Therefore 
 \[    r \hat{u} u_{s}t \otimes_{H} 1 =   r   u_{\sigma(s)} \tilde{u} t \otimes_{H} 1 = 
  r   u_{\sigma(s)}t  t^{-1}  \tilde{u} t \otimes_{H} 1 =    r   u_{\sigma(s)}t  \otimes_{H}  \phi( \tilde{u}) =   r   u_{\sigma(s)}t  \otimes_{H}  1 . \]
This means that 
 \[   r \hat{u} g \otimes_{H} \phi(h_{s})^{-1}   
 =  r  g  \otimes_{H}  \phi( h_{\sigma(s)})^{-1} ,  \]
and summing over $s$ shows that the action on the sum by $  r \hat{u}$ equals the action by $r$,  as expected.
 
 A candidate for a comultiplication $m^{*}$ on $R_{+}(G(a))$ has $(\alpha, a - \alpha)$-component of $m^{*}((H, \phi)^{G(a)})$ as, for each $t \in U_{\alpha, a - \alpha} \backslash G(a) / H$, the Lines
 $ g \otimes_{H} \phi(h)^{-1} $ where $g = uth$ ($u \in U_{\alpha, a - \alpha} , h \in H$) with $\phi = 1$ on 
 $t^{-1} U_{\alpha, a - \alpha} t \bigcap  H$.
 
 The previous calculation shows that, adding all these Lines over 
 \linebreak
 $t \in U_{\alpha, a - \alpha} \backslash G(a) / H$, $P_{\alpha, a - \alpha}/U_{\alpha, a - \alpha} \cong G(\alpha) \times G(a - \alpha)$, for each fixed double oset representative $t$, permutes these Lines - $\{  u_{s}t \otimes_{H} 1 \in  {\rm Ind}_{H}^{G(a)}(\phi) \} $ - so it is an element of $R_{+}(G(\alpha) \times G(a - \alpha))$.
 Write $\Sigma-{\rm orb}(t)$ for the set of permutation orbits of the action of $G(\alpha) \times G(a - \alpha)$ on 
 the set of Lines $\{  u_{s}t \otimes_{H} 1 \in  {\rm Ind}_{H}^{G(a)}(\phi) \} $ we receive, associated to the double coset of $t$ the sum of monomial objects
 \[  \oplus_{j \in \Sigma-{\rm orb}(t)}  ( G(j, \alpha, a - \alpha,t)  , 1)^{ G(\alpha) \times G(a - \alpha)} \]
 where $G(j, \alpha, a - \alpha,t)$ is the subgroup whose action gives the stabiliser of a chosen line in the 
 $ \Sigma-{\rm orb}(t)$ given by $j$.

  We have the formula
\[  m^{*}((H, \phi)^{G(a)})_{\alpha, a - \alpha}  = 
\Sigma_{t \in U_{\alpha, a - \alpha} \backslash G(a) / H} \   \oplus_{j \in \Sigma-{\rm orb}(t)}  \ ( G(j, \alpha, a - \alpha,t)  , 1)^{ G(\alpha) \times G(a - \alpha)} . \]
Also, in $R(G(\alpha) \times G(a - \alpha))$, $m^{*}b((H, \phi)^{G(a)}) $ is equal to the $U_{\alpha, a - \alpha} $-fixed subspace of $b(m^{*}((H, \phi)^{G(a)})_{\alpha, a - \alpha}  ))$.
\begin{question}{$_{}$}
\label{1.1a}
\begin{em}

Does $m^{*}$ extend to give a ``comultiplication'' (it does not map to the tensor product in the PSH manner but it does behave analogously) with respect to the multiplication $m$ on 
$R_{+ , \infty} = \oplus_{(a_{1}, \ldots , a_{n}) }  \ R_{+}(G(a_{1}) \times \ldots \times G(a_{n}))$?

It is clear that $m^{*}$ is associative.
\end{em}
\end{question} 

An example of the Hopf algebra conditions on $m^{*}$ and $m$ in $R_{1}$ is
$m^{*}m = (m \otimes m)(1 \otimes T \otimes 1)(m^{*} \otimes m^{*})$ so that
\[ m^{*}m : R(G(a)) \otimes R(G(b)) \longrightarrow R(G(a+b)) \longrightarrow R(G( \gamma + \delta)) \otimes R(G(a+b - \gamma - \delta)) \]
equals
\[ \begin{array}{c}
   R(G(a)) \otimes R(G(b))   \\
   \\
   \downarrow  \ m^{*} \otimes m^{*} \\
   \\
   R(G(\gamma)) \otimes R(G(a - \gamma)) \otimes  R(G(\delta)) \otimes R(G(b - \delta)) \\
   \\
      \downarrow  \ 1 \otimes T \otimes 1   \\
   \\
     R(G(\gamma)) \otimes  R(G(\delta))  \otimes R(G(a - \gamma))  \otimes R(G(b - \delta)) \\
   \\
      \downarrow  \ m \otimes m \\
   \\
   R(G( \gamma + \delta)) \otimes R(G(a+b - \gamma - \delta))  
\end{array} \]
when summed over $\gamma$ varying with each $\gamma + \delta$ fixed in the range
$1 \leq \gamma + \delta \leq a$. The $R_{+}(-)$-analogue is the equality of 
\[  R_{+}(G(a) \times G(b)) \stackrel{m}{\longrightarrow} R_{+}(G(a+b)) \stackrel{m^{*}}{\longrightarrow }
R_{+}(G( \gamma + \delta) \times G(a+b - \gamma - \delta)) \]
and
\[ \begin{array}{c}
R_{+}(G(a) \times G(b))  \\
\\ 
  \downarrow  \ m^{*}  \\
   \\
R_{+}(G(\gamma) \times G(a - \gamma) \times G(\delta) \times G(b - \delta)) \\
\\
     \downarrow  \ R_{+}(1 \times T \times 1)   \\
   \\
   R_{+}(G(\gamma) \times G(\delta)  \times G(a - \gamma)  \times G(b - \delta)) \\
\\
     \downarrow  \ m  \\
   \\
   R_{+}(G( \gamma + \delta) \times G(a+b - \gamma - \delta)) 
\end{array} \]
for fixed $\gamma + \delta$ and $\gamma$ summed over as above - at least when applied to 
$((H, \phi)^{G(a)}, (K, \psi)^{G(b)})$.

The $m^{*}m$-homomorphism sends $((H, \phi)^{G(a)}, (K, \psi)^{G(b)})$ first to 
\linebreak
$(H \times K, \phi \psi )^{G(a+b)}$ and then to 
\[ \Sigma_{t \in  \  U_{\gamma + \delta , a + b - \gamma - \delta}  \backslash G(a+b) / H \times K } \
 \oplus_{(j'') \in \Sigma-{\rm orb} } 
\  (G(j'', \gamma + \delta , a + b - \gamma - \delta,t) , 1)^{G(\gamma + \delta , a + b - \gamma - \delta) }. \]
The Lines in the general term in the above sum are $ g \otimes_{H} \phi(h)^{-1} $ where 
\linebreak
$g = uth \in G(a+b)$
for $h \in H, u \in   U_{\gamma + \delta , a + b - \gamma - \delta} $. Also the extra vanishing condition that $\phi \psi = 1$ on $t^{-1}  U_{\gamma + \delta , a + b - \gamma - \delta}  t \bigcap H \times K$ holds if and only if 
$\phi  = 1$ on $t^{-1}  U_{\gamma + \delta , a + b - \gamma - \delta}  t \bigcap H $ and $\psi  = 1$ on $t^{-1}  U_{\gamma + \delta , a + b - \gamma - \delta}  t \bigcap  K$.

The other map sends $((H, \phi)^{G(a)}, (K, \psi)^{G(b)})$ to 
\[ \Sigma_{t', t''} \  \oplus_{(j.j') \in \Sigma-{\rm orb} \ {\rm pairs}}  \ ((G(j, \gamma, a - \gamma, t'), 1)^{G(j', \gamma, a - \gamma)} , (G(\delta, b - \delta, t'') ,1)^{G(\delta, b - \delta)} ) \]
where $t' \in U_{\gamma, a - \gamma} \backslash \ G(a)/ H$ and $t'' \in U_{\delta, b - \delta} \backslash \ G(b)/ K$ which is then mapped to 
\[ \Sigma_{t', t''} \ 
 \oplus_{(j'') \in \Sigma-{\rm orb} } \ (( G(j'', \gamma + \delta, a + b -  \gamma - \delta), t' t'') , 1 )^{G(\gamma + \delta, a + b -  \gamma - \delta)} ) . \]
However, fixing $\gamma + \delta$ and summing over varying $\gamma$, the combinatorics of double cosets in the general linear groups of a field, given in (\cite{AVZ81} Appendix A3.5 p.173) show that the above two compositions send $((H, \phi)^{G(a)}, (K, \psi)^{G(b)})$  to the same element, as expected.

In general similar double coset combinatorics will yield the commutativity of any plausible diagram involving compositions of $m^{*}$'s and $m$'s. Typically
\newline
\begin{picture}(150,100)(-60,0)
\put(125,90){ $ m $ }
\put(100,20){ $ m \cdot {\rm shuffle} $ }
\put(60,55){$ m^{*} $}
\put(180,55){$ m^{*} $}
\put(-40,80){$ R_{+}(G(a_{1})) \otimes \ldots \otimes  R_{+}(G(a_{n}))$}
\put(55,70){\vector(0,-1){30}}
\put(-30,10){$ \Sigma_{\beta} \ R_{+}(G(\beta_{1}) \times \ldots )  $}
\put(100,12){\vector(1,0){25}}
\put(120,82){\vector(1,0){25}}
\put(150,80){$R_{+}(G(a_{1} + \ldots ) \times \ldots \ \times G(\ldots + a_{n}))  $}
\put(200,70){\vector(0,-1){30}}
\put(160,10){$ R_{+}(G(\gamma) \times G(a_{1}+ \ldots - \gamma) \times \ldots )$}
\end{picture}
\newline
where $\beta$ runs through all ordered sub-partitions of the those appearing in the lower right group.

We have another batch of  commutative diagrams of the following typical form.
\newline
\begin{picture}(150,100)(-60,0)
\put(125,85){ $ b $ }
\put(130,20){ $ b(-)^{U} $ }
\put(60,55){$m^{*}$}
\put(180,55){$ m^{*} $}
\put(-10,80){$   R_{+}(G(a_{1}) \times \ldots )  $}
\put(55,70){\vector(0,-1){30}}
\put(-30,10){$R_{+}(G(\gamma) \times G(a_{1} - \gamma) \times \ldots ) $}
\put(130,12){\vector(1,0){25}}
\put(120,82){\vector(1,0){25}}
\put(180,80){$  R(G(a_{1}) \times \ldots ) $}
\put(200,70){\vector(0,-1){30}}
\put(180,10){$ R(G(\gamma) \times G(a_{1} - \gamma) \times \ldots ) $}
\end{picture}
\newline
where $ b(-)^{U} $ denotes $b$ followed by taking $U(\gamma, a - \gamma)$-fixed vectors.

\underline{{\bf Orbit multiplicities}}

For $t \in U_{\alpha, a - \alpha} \backslash G(a) / H$ we have the set of $ \Sigma-{\rm orb}(t)$ which divide up the set of spacial Lines associated to $t$ and the multiplicity of a given orbit $j$, denoted by $\mu(j)$ will be the number of orbits conjugate  to $j$ in this set so that another form for the formula is
\[  \begin{array}{l}
m^{*}((H, \phi)^{G(a)})_{\alpha, a - \alpha} \\
\\
 = 
\Sigma_{t \in U_{\alpha, a - \alpha} \backslash G(a) / H} \   \oplus_{ \stackrel{j \in \Sigma-{\rm orb}(t)}{  {\rm mod \ conjugation}} } \ \mu(j) ( G(j, \alpha, a - \alpha,t)  , 1)^{ G(\alpha) \times G(a - \alpha)} .
\end{array}   \]

\underline{{\bf The hyperHecke relation in ${\mathcal H}_{cmc}(G(\alpha) \times G(a - \alpha))$}}

Consider the relation involving 
\[  [G(j', \alpha , a - \alpha, t) , 1) , g , G(j'', \alpha , a - \alpha, tg^{-1)} , 1))  \in   {\mathcal H}_{cmc}(G(\alpha) \times G(a - \alpha)),  \]
which is given in the following ection on the definition of the hyperHecke algebra. This relation respects our construction above. This is because, for example,  
\[  [G(j', \alpha , a - \alpha, t) , 1) , gk , G(j'', \alpha , a - \alpha, tg^{-1)} , 1)) \]
with $k \in G(j', \alpha , a - \alpha, t) , 1) $ then $``\psi(k) =1$'' because $``\psi =1$ on $G(j', \alpha , a - \alpha, t) , 1) $ in the current case where $\underline{\phi} =1$.

{\bf Alternatively:} \ I could have replaced $\phi$ by ${\rm Hom}_{{\mathbb F}_{q}^{*}}(\underline{\phi}, \phi)$  which has a trivial action by the central subgroup of scalars.

\section{A coproduct for the hyperHecke algebra of $GL_{*}{\mathbb F}_{q}$?}

\underline{{\bf The hyperHecke algebra }}

The following is recalled from (\cite{Sn20} \S2). In my treatment of local fields and admissible representations I have usually treated all representations with a fixed central character $\underline{\phi}$.  This central character prevents the ``$\phi=1 $ on $tU_{\alpha, a - \alpha}t^{-1} \bigcap g^{-1}Hg$'' if it is non-trivial. Therefore, pro tem I shall either ignore the central character of pretend it is trivial.
 
 Let $G$ be a locally profinite group and let $k$ be an algebraically closed field. Suppose that $\underline{\phi} : Z(G) \longrightarrow k^{*}$ is a fixed $k$-valued, continuous character on the centre $Z(G)$ of $G$. Let ${\mathcal M}_{cmc, \underline{\phi}}(G)$ be the poset of pairs $(H, \phi)$ where $H$ is a subgroup of $G$, containing $Z(G)$, which is compact, open modulo the centre of $G$ and $\phi : H \longrightarrow k^{*}$ is a $k$-valued, continuous character whose restriction to $Z(G)$ is $\underline{\phi}$\footnote{As mentioned a few lines earlier, $\underline{\phi}$ is temporarily assumed trivial.}.
 
 We define the hyperHecke algebra, ${\mathcal H}_{cmc}(G)$,  to be the $k$-algebra given by the following generators and relations. For $(H, \phi), (K, \psi) \in {\mathcal M}_{cmc,\underline{\phi}}(G)$,  write 
$[(K, \psi), g, (H, \phi)]$ for any triple consisting of $g \in G$, characters $\phi, \psi$ on 
subgroups $H, K \leq G$, respectively such that 
\[    (K, \psi) \leq (g^{-1}Hg, (g)^{*}(\phi)) \]
which means that $K \leq  g^{-1}Hg$ and that $\psi(k) = \phi(h)$ where $k = g^{-1}hg$ for
 $h \in H, k \in K$.

 Let ${\mathcal H}$ denote the $k$-vector space with basis given by these triples. Define a product on these triples by the formula 
\[  [(H, \phi), g_{1}, (J, \mu)]  \cdot  [(K, \psi), g_{2}, (H, \phi)] =   [(K, \psi), g_{1}g_{2}, (J, \mu)]  \]
and zero otherwise. This product makes sense because 

(i)  \   if $K \leq g_{2}^{-1} H g_{2}$ and
 $H \leq g_{1}^{-1} J g_{1}$ then $K \leq  g_{2}^{-1} H g_{2} \leq   g_{2}^{-1} g_{1}^{-1} J g_{1} g_{2} $ 
 
 and 
 
 (ii)  \  if $\psi(k) = \phi(h) = \mu(j), $ where $k = g_{2}^{-1}hg_{2},  h = g_{1}^{-1}j g_{1}$ then
 \linebreak
  $k = g_{2}^{-1}  g_{1}^{-1}j g_{1}   g_{2}$. 
 
 This product is clearly associative and we define an algebra ${\mathcal H}_{cmc}(G)$ to be ${\mathcal H}$ modulo the relations ( \cite{Sn18}\footnote{For the purposes of this essay I  am using throughout a convention in which $g$ is replaced by $g^{-1}$ to make the composition coincide with the conventions of induced representations which I shall use here.})
 \[   [(K, \psi), gk, (H, \phi)]  = \psi(k^{-1}) [(K, \psi), g, (H, \phi)]  \]
and
 \[     [(K, \psi), hg, (H, \phi)]  = \phi(h^{-1}) [(K, \psi), g, (H, \phi)] .    \]

Now let us consider the finite case where $G = GL_{a}{\mathbb F}_{q}$ denoted by $G(a)$. For $(H, \phi), (K, \psi) \in {\mathcal M}_{cmc,\underline{\phi}}(G(a))$,  recall that
$[(K, \psi), g, (H, \phi)]$ denotes a triple consisting of $g \in G$, characters $\phi, \psi$ on 
subgroups $H, K \leq G(a)$, respectively such that 
\[    (K, \psi) \leq (g^{-1}Hg, (g)^{*}(\phi)) \]
which means that $K \leq  g^{-1}Hg$ and that $\psi(k) = \phi(h)$ where $k = g^{-1}hg$ for
 $h \in H, k \in K$.

 We have $G(\alpha) \times G(a - \alpha) \leq G(a)$. From the previous section  we have the formula
\[  \begin{array}{l}
m^{*}((H, \phi)^{G(a)})_{\alpha, a - \alpha} \\
\\
 = 
\Sigma_{t \in U_{\alpha, a - \alpha} \backslash G(a) / H} \   \oplus_{ \stackrel{j \in \Sigma-{\rm orb}(t)}{  {\rm mod \ conjugation}} } \ \mu(j) ( G(j, \alpha, a - \alpha,t)  , 1)^{ G(\alpha) \times G(a - \alpha)} .
\end{array}   \]

Also, in $R(G(\alpha) \times G(a - \alpha))$, $m^{*}b((H, \phi)^{G(a)}) $ is equal to the $U_{\alpha, a - \alpha} $-fixed subspace of $b(m^{*}((H, \phi)^{G(a)})_{\alpha, a - \alpha}  ))$.
In this formula, as an induced representation, $ (G(\alpha, a - \alpha,t)  , 1)^{ G(\alpha) \times G(a - \alpha)}$
consists of the Lines $ g' \otimes_{H} \phi(h)^{-1} $ where $g' = uth$ ($u \in U_{\alpha, a - \alpha} , h \in H$) with $\phi = 1$ on $t U_{\alpha, a - \alpha} t^{-1} \bigcap  H$.

The condition that $g^{*}(\phi)=1$ on $tU_{\alpha, a - \alpha}t^{-1} \bigcap g^{-1}Hg$ is equivalent to the condition that $\phi = 1$ on $gt U_{\alpha, a - \alpha}t^{-1}g^{-1}  \bigcap H$ and, in addition, it implies that $\psi=1$ on $tU_{\alpha, a - \alpha}t^{-1} \bigcap K$.

We have a well-defined left  $k[G]$-module homomorphism (\cite{Sn20}  Appendix \S12)
\[    [(K, \psi), g, (H, \phi)] : k[G] \otimes_{k[K]} k_{\psi}  \longrightarrow  k[G] \otimes_{k[H]} k_{\phi}  \]
given by the formula $[(K, \psi), g, (H, \phi)](g' \otimes_{k[K]} v) = g' g^{-1} \otimes_{k[H]} v$. 

By comparing the ``comultiplication'' on the objects in the previous section we shall define a `comultiplication'' on the hyperHecke algebra of $G(a)$ whose $G(\alpha) \times G(a - \alpha)$-component will 
satisfy $ m^{*} [(K, \psi), g, (H, \phi)] = 0  $ except, for $t \in U_{\alpha, a - \alpha} \backslash G/ K$, when $\phi = 1$ on $gt U_{\alpha, a - \alpha}t^{-1}g^{-1}  \bigcap H$ in which case we shall proceed as follows. In this case $\psi =1$ on  $t U_{\alpha, a - \alpha}t^{-1}g  \bigcap K$ and $  [(K, \psi), g, (H, \phi)] $ gives a 
$G(\alpha) \times G(a - \alpha)$-permutation module homomorphism from the direct sum of the permutation modules in the set 
$j' \in \Sigma-{\rm orb}(t)$ to the sum in the set $j'' \in \Sigma-{\rm orb}(tg^{-1})$. This map will give the sum in the comultiplication formula
\[  \begin{array}{l}
m^{*}_{G(\alpha) \times G(a - \alpha)} [(K, \psi), g, (H, \phi)] \\
\\
= \sum_{j', j''} \ [ ( G(j', \alpha, a - \alpha,t)  , 1) , g , ( G(j'', \alpha, a - \alpha,tg^{-1})  , 1) ] )  . 
\end{array} \]

\section{A product for the hyperHecke algebra of $GL_{*}{\mathbb F}_{q}$?}

Persisting with the assumption that the central character $\underline{\phi}$ is trivial, we define a ``parabolic induction product''
\[ \begin{array}{c}
 {\mathcal H}_{cmc}(G(a_{1}))  \otimes  \ldots \otimes {\mathcal H}_{cmc}(G(a_{r})) \otimes 
 {\mathcal H}_{cmc}(G(b_{1}))  \otimes  \ldots \otimes {\mathcal H}_{cmc}(G(b_{r}))  \\
 \\
 \downarrow  \hspace{5pt} m  \\
 \\
  {\mathcal H}_{cmc}(G(a_{1}+ b_{1}))  \otimes  \ldots \otimes {\mathcal H}_{cmc}(G(a_{r}+ b_{r}))
\end{array} \]
by the formula
\[ \begin{array}{l}
m(  [(K_{1}, \psi_{1}), g_{1}, (H_{1}, \phi_{1})] \otimes \ldots \otimes  [(K'_{1}, \psi'_{1}), g'_{1}, (H'_{1}, \phi'_{1})] \otimes \ldots) \\
\\
=  [(K_{1} \times K'_{1}, ( \psi_{1} \cdot \psi'_{1})), (g_{1}, g'_{1}) ,(H_{1} \times H'_{1}, (\phi_{1} \cdot \phi'_{1})) ] \otimes \ldots .
\end{array} \]
The multiplication on
\[ \begin{array}{l}
 {\mathcal H}_{cmc}(G(a_{1}))  \otimes  \ldots \otimes {\mathcal H}_{cmc}(G(a_{r})) \otimes 
 {\mathcal H}_{cmc}(G(b_{1}))  \otimes  \ldots \otimes {\mathcal H}_{cmc}(G(b_{s}))  
 \end{array} \]
is defined to vanish when $r \not= s$.

\underline{{\bf Non-trivial central character}}

In this case we simply intersect $U_{\alpha, a - \alpha}$ with $SL_{a}{\mathbb F}_{q} \leq G(a)$. That is, $(-)^{U_{\alpha, a - \alpha}}$ is replaced by $(-)^{U_{\alpha, a - \alpha} \bigcap SL_{a}{\mathbb F}_{q}}$
and the condition that $\phi = 1$ on $g^{-1}t U_{\alpha, a - \alpha} tg^{-1} \bigcap H$ becomes 
$\phi = 1$ on $g^{-1}t (U_{\alpha, a - \alpha} tg^{-1}  \bigcap SL_{a}{\mathbb F}_{q} \bigcap H$.

Here I have simply said ``intersect with $SL$'' meaning, for larger products of $GL$'s, to intersect with corresponding products of $SL_{a_{i}}$'s.

In the case of local fields each $H$ is compact open modulo the centre (which is the scalar matrices). Hence $H$ is in a product ${\rm scalars} \times H'$ where $H'$ is a compact open subgroup. We can construct the hyperHecke coproduct using the $\phi = 1$ condition on $g^{-1}t (U_{\alpha, a - \alpha} tg^{-1}  \bigcap SL_{a}{\mathbb F}_{q} \bigcap H'$ and then extend the hyperHecke coproduct triples to ${\rm scalars} \times H'$ 
by means of which $U_{\alpha, a - \alpha}$ acts trivially and the scalars act via $\underline{\phi}$. This induces the required $G(\alpha) \times G(a - \alpha)$-action with $\underline{\phi}$-central character action. Hence we obtain a coproduct on the local $GL$'s hyperHecke algebra.

In the case of finite fields the analogous construction is clear and works to give a hyperHecke ``coproduct''

In the composition
\[ \begin{array}{c}
 {\mathcal H}_{cmc}(G(a_{1}))   \otimes 
 {\mathcal H}_{cmc}(G(b_{1}))    \\
 \\
 \downarrow  \hspace{5pt} m  \\
 \\
  {\mathcal H}_{cmc}(G(a_{1}+ b_{1}))  \\
  \\
   \downarrow  \hspace{5pt} m^{*}  \\
 \\
 {\mathcal H}_{cmc}(G(\alpha))   \times 
 {\mathcal H}_{cmc}(G(a_{1} + b_{1} - \alpha ))    \\
\end{array} \]
mapping the element $[(K_{1}, \psi_{1}), g_{1}, (H_{1}, \phi_{1})] \otimes \ldots \otimes  [(K'_{1}, \psi'_{1}), g'_{1}, (H'_{1}, \phi'_{1})] $ the condition that $\phi_{1} \cdot \phi'_{1} = 1$ holds if and only if it holds separately for $\phi_{1} $ and $ \phi'_{1}$. This fact, together with the usual $GL$-double coset combinatorics, yields the following commutative diagram 
\newline
\begin{picture}(150,100)(-60,0)
\put(125,90){ $ m $ }
\put(100,-15){ $ \Sigma_{\delta_{1}, \delta_{2}} \ (- \cdot - )  $ }
\put(60,55){$ m^{*} \otimes m^{*} $}
\put(180,55){$ m^{*} $}
\put(-40,80){$  {\mathcal H}_{cmc}(G(a_{1}))   \otimes 
 {\mathcal H}_{cmc}(G(b_{1}))  $}
\put(55,70){\vector(0,-1){30}}
\put(-50,25){$ {\mathcal H}_{cmc}(G(\delta_{1} ))   \times 
 {\mathcal H}_{cmc}(G(a_{1} + b_{1} - \delta_{1} )) $}
\put(-50,0){$ \otimes {\mathcal H}_{cmc}(G(\delta_{2} ))   \times 
 {\mathcal H}_{cmc}(G(a_{1} + b_{1} - \delta_{2} )) $}
\put(140,12){\vector(1,0){25}}
\put(120,82){\vector(1,0){25}}
\put(150,80){$  {\mathcal H}_{cmc}(G(a_{1}+ b_{1}))   $}
\put(200,70){\vector(0,-1){30}}
\put(170,15){$\Sigma_{\delta_{1}, \delta_{2} } \  {\mathcal H}_{cmc}(G(*))   
 $}
 \put(210,0){$    \times 
 {\mathcal H}_{cmc}(G(a_{1} + b_{1} - * )) $}
\end{picture}
\newline
\vspace{10pt} 

when summed over all $\delta_{1}, \delta_{2}$ such that $\alpha = \delta_{1}+ \delta_{2}$. 

Similar commutative PSH-like $m^{*}$, $m$ diagrams hold for longer ordered partitions.

\section{Relation with the hyperHecke algebra composition product}

 In $ {\mathcal H}_{cmc}(G(a))  $ consider applying the coproduct $m^{*}$ to the composition product
\[  [(H, \phi), g_{1}, (J, \mu)]  \cdot  [(K, \psi), g_{2}, (H, \phi)] =   [(K, \psi), g_{1}g_{2}, (J, \mu)] . \]

If $\phi \mu = 1$ on $(g_{1}g_{2})t U_{\alpha, a - \alpha} t (g_{1}g_{2})^{-1} \bigcap SL \bigcap J$ then $\mu = 1$ on 
\linebreak
$g_{1}t U_{\alpha, a - \alpha} t g_{1}^{-1} \bigcap SL \bigcap H$ and $  [(K, \psi), g_{1}g_{2}, (J, \mu)] $ is the $G(\alpha ) \times G(a - \alpha)$- composition of maps of the permutation modules given by the special Lines - so the factors in the composition , written as ${\mathbb Z}$-linear combinations of hyperHecke algebra elements, is the composition product of the matrices of hyperHecke algebra elements. Hence
\[ \begin{array}{l}
m^{*}( [(H, \phi), g_{1}, (J, \mu)]  \cdot  [(K, \psi), g_{2}, (H, \phi)] ) \\
\\
=  m^{*}( [(H, \phi), g_{1}, (J, \mu)]  ) \cdot  m^{*}(   [(K, \psi), g_{2}, (H, \phi)] ) . 
\end{array} \]
In other cases bth sides agree because they are zero.

So $m^{*}$ is a homomorphism of hyperHecke algebra. We have a the following commutative diagram.
\newline
\begin{picture}(150,100)(-60,0)
\put(105,90){ $ (- \cdot - ) $ }
\put(100,20){ $  (- \cdot - )  $ }
\put(60,55){$ m^{*} \otimes m^{*} $}
\put(210,55){$ m^{*}_{\alpha, a - \alpha}$}
\put(-40,80){$  {\mathcal H}_{cmc}(G(a))   \otimes 
 {\mathcal H}_{cmc}(G(a)) $}
\put(55,70){\vector(0,-1){30}}
\put(-50,25){$ {\mathcal H}_{cmc}(G(\alpha ))   \times 
 G(a  - \alpha )) $}
\put(-50,0){$ \otimes{\mathcal H}_{cmc}(G(\alpha ))   \times 
 G(a  - \alpha ))  $}
\put(100,12){\vector(1,0){25}}
\put(120,82){\vector(1,0){25}}
\put(170,80){$  {\mathcal H}_{cmc}(G(a))   $}
\put(200,70){\vector(0,-1){30}}
\put(160,10){$ {\mathcal H}_{cmc}(G(\alpha ))   \times 
 G(a  - \alpha ))  $}
\end{picture}
\newline

\underline{{\bf Multiplication and composition product}}

Consider the diagram
\newline
\begin{picture}(150,100)(-60,0)
\put(90,95){ $ (m \otimes m ){\rm shuffle} $ }
\put(100,20){ $  (- \cdot - )  $ }
\put(210,55){$ (-,-)$}
\put(-40,80){$  {\mathcal H}_{cmc}(G(a))   \otimes 
 {\mathcal H}_{cmc}(G(a)) $}
 \put(-40,60){$   \otimes  {\mathcal H}_{cmc}(G(b))   \otimes 
 {\mathcal H}_{cmc}(G(b)) $}
 \put(-25,40){$(-,-) \otimes (-,-)$}
\put(65,50){\vector(0,-1){15}}
\put(-50,10){$ {\mathcal H}_{cmc}(G(a ))   \otimes  {\mathcal H}_{cmc}(G(b ))
 $}
\put(100,12){\vector(1,0){25} }
\put(120,82){\vector(1,0){25}}
\put(150,80){$  {\mathcal H}_{cmc}(G(a+b)) \otimes  {\mathcal H}_{cmc}(G(a+b))   $}
\put(200,70){\vector(0,-1){30}}
\put(160,10){$ {\mathcal H}_{cmc}(G(\alpha )   \times 
 G(b ))  $}
\end{picture}
\newline

Going round by the counter-clockwise route we have
\[ \begin{array}{l}
[(H, \phi), g_{1}, (J, \mu)] \otimes  [(K, \psi), g_{2}, (H, \phi)] \otimes  \\
\\
\hspace{50pt} [ (H, \phi_{1}) , g_{3}, (J, \mu_{1}] \otimes 
 [(K_{1}, \psi_{1}), g_{4}, (H_{1}, \phi_{1})] \\
 \\
 \mapsto  [(K, \psi), g_{1}g_{2}, (J, \mu)] \otimes  [(K_{1}, \psi_{1}), g_{3}g_{4}, (J_{1}, \mu_{1})]   \\
 \\
 \mapsto  [(K \times K_{1} , \psi \psi_{1}) , (g_{1}g_{3}, g_{2}g_{4}), (J \times J_{1} , \mu \mu_{1})] 
\end{array} \]
which is the outcome of going round by the clockwise route, making the diagram commutative.

\section{Boltje's bilinear form \cite{Bo90} }

Boltje defined a non-symmetrical bilinear form (\cite{Sn94} p.46) 
\[ [ - , -, ] : R_{+}(G) \times R_{+}(G) \longrightarrow {\mathbb Z} \]
by the formula (\cite{Sn94} p.34) 
\[  [ (H, \phi)^{G} , (H', \phi')^{G}   ] = \# \{ g \in H \backslash G / H' \ | \ (H, \phi) \leq (gH'g^{-1}, (g^{-1})^{*}( \phi') \} .  \]
It satisfies the adjunction relation
\[  [ x , a_{G}(\rho) ] = \langle b_{G}(x), \rho \rangle \]
where $a_{G}(-) $ is the explicit Brauer induction map of ( \cite{Bo90} ; see also \cite{Sn94}).

Thus $[-,-]$ is a non-symmetric bilinear form on objects in $R_{+}(G)$. On objects in $R(G)$ we have the Schur inner product.

\section{$c-{\rm Ind}_{H}^{GL_{n}F}(k_{\phi})$ with $F$ a $p$-adic local field of characteristic zero}  

Let $G$ be the locally profinite group $GL_{n}F$. 
Let us begin by recalling, from (\cite{Sn18} Chapter Two \S1)\footnote{In fact, I am using the description and conventions of \cite{Sn20} \S5 concerning induced and compactly induced smooth representations.}, induced and compactly induced smooth representations. 
\begin{definition}{$_{}$}
\label{5.1A}
\begin{em}

Let $H \subseteq G$ a closed subgroup. Thus $H$ is also  locally profinite. Let 
\[   \phi : H \longrightarrow  k^{*}     \]
be a smooth representation of $H$ which extends $\underline{\phi}$ on $F^{*} \leq H$.  Set $X$ equal to the space of functions $f: G \longrightarrow k_{\phi}  $ such that (writing simply $h \cdot w$ for $\phi(h)(w)$ if $h \in H, w \in  k_{\phi} $)

(i) \  $f(hg) = h \cdot f(g)$ for all $h \in H, g \in G$,

(ii)  \  there is a compact open subgroup $K_{f} \subseteq G$ such that $f(gk) = f(g)$ for all $g \in G, k \in K_{f}$.

The (left) action of $G$ on $X$ is given by $(g \cdot f)(x)= f(xg)$ and
\[   \Sigma :  G  \longrightarrow  {\rm Aut}_{k}(X)  \]
gives a smooth representation of $G$ (in the classical smoothness where the residue characteristic of $F$ and $k$ are different).

The representation $\Sigma$ is called the representation of $G$ smoothly induced from $\sigma$ and is usually denoted by $\Sigma = {\rm Ind}_{H}^{G}(k_{\phi})$.
\end{em}
\end{definition}

\begin{dummy}
\label{5.2A}
\begin{em}

\[  (g \cdot f)(hg_{1}) = f(hg_{1}g) = h f(g_{1}g) = h (g \cdot f)(g_{1}) \]  
so that $(g \cdot f)$ satisfies condition (i) of Definition \ref{5.1A}. 

Also, by the same discussion as in the finite group case (see Appndix \S4), the formula will give a left $G$-representation, providing that $g \cdot f \in X$ when $f \in X$.
However, condition (ii) asserts that there exists a compact open subgroup $K_{f}$ such that $k \cdot f = f$ for all $k \in K_{f}$. The subgroup $gK_{f}g^{-1}$ is also a compact open subgroup and, if $k \in K_{f}$,
we have
\[ (gkg^{-1}) \cdot (g \cdot f) = (gkg^{-1}g) \cdot f =  (gk) \cdot f = (g \cdot (k \cdot f)) =  (g \cdot f) \]
so that $g \cdot f \in X$, as required.

The smooth representations of $G$ form an abelian category ${\rm Rep}(G)$.
\end{em}
\end{dummy}

\begin{proposition}{$_{}$}
\label{5.3A}
\begin{em}

The functor
\[  {\rm Ind}_{H}^{G} : {\rm Rep}(H) \longrightarrow   {\rm Rep}(G)  \]
is additive and exact.
\end{em}
\end{proposition}
\begin{proposition}{(Frobenius Reciprocity) }
\label{5.4A}
\begin{em}

There is an isomorphism
\[  {\rm Hom}_{G}( \pi,  {\rm Ind}_{H}^{G}(\sigma)) \stackrel{\cong}{\longrightarrow}  {\rm Hom}_{H}( \pi,  \sigma)  \]
given by $\phi \mapsto \alpha \cdot \phi$ where $\alpha$ is the $H$-map  
\[   {\rm Ind}_{H}^{G}(\sigma)  \longrightarrow  \sigma  \]
given by $\alpha(f) = f(1)$.
\end{em}
\end{proposition} 

\begin{dummy}
\label{5.5A}
\begin{em}

In general, if $H \subseteq Q$ are two closed subgroups there is a $Q$-map
\[   {\rm Ind}_{H}^{G}(\sigma)  \longrightarrow   {\rm Ind}_{H}^{Q}(\sigma)   \]
given by restriction of functions. Note that $\alpha$ in Proposition \ref{5.4A} is the special case where $H=Q$.
\end{em}
\end{dummy}

\begin{dummy}{The c-Ind variation (\cite{Sn20} \S 5.6)}
\label{5.6A}
\begin{em}

Inside $X$ let $X_{c}$ denote the set of functions which are compactly supported modulo $H$. This means that the image of the support
\[   {\rm supp}(f) = \{ g \in G \ | \  f(g) \not= 0  \}  \]
has compact image in $H \backslash G$. Alternatively  there is a compact subset $C \subseteq G$ such that $  {\rm supp}(f) \subseteq H \cdot C$.

The $\Sigma$-action on $X$ preserves $X_{c}$, since $ {\rm supp}(g \cdot f) =  {\rm supp}(f) g^{-1} \subseteq HCg^{-1}$, and we obtain $X_{c} =  c- {\rm Ind}_{H}^{G}(W) $, the compact induction of $W$ from $H$ to $G$.

This construction is of particular interest when $H$ is open. There is a canonical left $H$-map (see the Appendix  in induction in the case of finite groups)
\[    f : W \longrightarrow  c- {\rm Ind}_{H}^{G}(W)  \]
given by $w \mapsto f_{w}$ where $f_{w}$ is supported in $H$ and $f_{w}(h) = h \cdot w$ (so $f_{w}(g) =0$ if $g \not\in H$).

For $g \in G$ we have
\[ \begin{array}{ll}
(g \cdot f_{w})(x) = f_{w}(xg) & =  \left\{  \begin{array}{cc}
0 & {\rm if } \ xg \not\in H, \\
\\
(xg^{-1}) \cdot w & {\rm if } \ xg \in H, \\
\end{array} \right.   \\
\\
& =  \left\{  \begin{array}{cc}
0 & {\rm if } \  x  \not\in Hg^{-1}, \\
\\
(xg^{-1}) \cdot w & {\rm if } \ x \in Hg^{-1}. \\
\end{array} \right. 
\end{array} \]

We shall be particularly interested in the case when ${\rm dim}_{k}(W) =1$. In this case we write $W = k_{\phi}$ where $\phi : H \longrightarrow k^{*}$ is a continuous/smooth character and, as a vector space with a left $H$-action $W =k$ on which $h \in H$ acts by multiplication by $\phi(h)$. In this case $\alpha_{c}$ is an
injective left $k[H]$-module homomorphism of the form
\[   f : k_{\phi} \longrightarrow   c- {\rm Ind}_{H}^{G}(k_{\phi}) .\]

\end{em}
\end{dummy}

\begin{lemma}{(\cite{Sn20} \S 5.7)}
\label{5.7A}
\begin{em}

Let $H$ be an open subgroup of $G$. Then

(i) \  $f : w \mapsto f_{w}$ is an $H$-isomorphism onto the space of functions $f \in c- {\rm Ind}_{H}^{G}(W) $ such that $  {\rm supp}(f)  \subseteq H$.

(ii) \  If $w \in W$ and $h \in H$ then $   h \cdot f_{w}   =  f_{h \cdot w}$.

(iii)  \   If ${\mathcal W}$ is a $k$-basis of $W$ and ${\mathcal G}$ is a set of coset representatives for 
$H \backslash G  $ then
\[   \{  g  \cdot f_{w}  \  |  \  w \in {\mathcal W}, \   g \in {\mathcal G}  \}\]
is a $k$-basis of $c- {\rm Ind}_{H}^{G}(W) $.
\end{em}
\end{lemma}
\vspace{2pt}

\begin{example}  (\cite{Sn20} \S 5.8)
\label{5.8}
\begin{em}

Let $K$ be a $p$-adic local field with valuation ring ${\mathcal O}_{K}$ and $\pi_{K}$ a generator of the maximal ideal of ${\mathcal O}_{K}$. Suppose that $G = GL_{n}K$ and that $H$ is a subgroup containing the centre of $G$ (that is, the scalar matrices $K^{*}$).  If $H$ is compact, open modulo $K^{*}$ then there is a subgroup $H'$ of finite index in $H$ such that $H' = K^{*}H_{1}$ with $H_{1}$ compact, open in $SL_{n}K$. This can be established by studying the simplicial action of $GL_{n}K$ on a suitable barycentric subdivision of the Bruhat-Tits building of $SL_{n}K$ (see \cite{Sn18} Chapter Four \S1).

To show that $H$ is both open and closed it suffices to verify this for $H'$. Firstly $H'$ is open, since it is $H' = \bigcup_{z \in K^{*}} \ zH_{1} = \bigcup_{s \in {\mathbb Z}}  \  \pi_{K}^{s} H_{1}$. 

Also $H' = K^{*}H_{1}$ is closed. Suppose that $X'  \not\in   K^{*}H_{1}$. $K^{*}H_{1}$ is closed under mutiplication by the multiplicative group generated by $\pi_{K}$ so that $\pi_{K}^{m} X' \not\in  K^{*}H_{1}$ for all $m$. By conjugation we may assume that $H_{1}$ is a subgroup of $SL_{n}{\mathcal O}_{K}$, which is the maximal compact open subgroup of $SL_{n}K$, unique up to conjugacy. Choose the smallest non-negative integer $m$ such that every entry of $X = \pi_{K}^{m} X' $ lies in
$ {\mathcal O}_{K}$. Therefore we may write $0 \not= {\rm det}(X) = \pi_{K}^{s} u$ where $u \in {\mathcal O}_{K}^{*}$ and $1 \leq s$. Now suppose that $V$ is an $n \times  n$ matrix with entries in $ {\mathcal O}_{K}$ such that $X + \pi_{K}^{t}V  \in K^{*}H_{1}$. Then
\[ {\rm det}(X + \pi_{K}^{t}V)  \equiv  \pi_{K}^{s} u  \ ({\rm modulo} \  \pi_{K}^{t}) . \]
So that if $t > s$ then $s$ must have the form $s = nw$ for some integer $w$ and 
$\pi_{K}^{-w} (X + \pi_{K}^{t}V) \in GL_{n}{\mathcal O}_{K} \bigcap K^{*}H_{1} = H_{1}$. Therefore all the entries in $\pi_{K}^{-w} X$ lie in ${\mathcal O}_{K} $ and $\pi_{K}^{-w} X \in GL_{n}{\mathcal O}_{K}$. Enlarging $t$, if necessary, we can ensure that $\pi_{K}^{-w} X \in H_{1}$, since $H_{1}$ is closed (being compact), and therefore 
$X' \in K^{*} H_{1}$, which is a contradiction.

Since $H $ is both closed and open in $GL_{n}K$ we may form the admissible representation $c-{\rm Ind}_{H}^{GL_{n}K}(k_{\phi})$ for any continuous character $\phi : H \longrightarrow k^{*}$ and apply Lemma \ref{5.7A}. 

If $g \in GL_{n}K, h \in H$ then $(g \cdot f_{1})(x) =  \phi(xg)$ if $xg \in H$ and zero otherwise. On the other hand, $(gh \cdot f_{1})(x) =  \phi(xgh) = \phi(h)\phi(xg)$ if $xg \in H$ and zero otherwise. Therefore as a left $GL_{n}K$-representation $c-{\rm Ind}_{H}^{GL_{n}K}(k_{\phi})$
is isomorphic to 
\[ k[GL_{n}K] /( \phi(h) g  -  gh \ | \  g \in GL_{n}K, \ h \in H)  \]
with left action induced by $g_{1} \cdot g = g_{1}g$. 

This vector space is isomorphic to the $k$-vector space whose basis is given by $k$-bilinear tensors over $H$ of the form $g \otimes_{k[H]} 1$ as in the case of finite groups. The basis vector $g \cdot f_{1}$ corresponds to $g \otimes_{H} 1$ and $GL_{n}K$ acts on the tensors by left multiplication, as usual (see Appendix \S4 in the finite group case).
\end{em}
\end{example}

\begin{proposition}{$_{}$}
\label{5.9}
\begin{em}

The functor
\[   c- {\rm Ind}_{H}^{G} : {\rm Rep}(H) \longrightarrow   {\rm Rep}(G)  \]
is additive and exact.
\end{em}
\end{proposition}

\begin{proposition}{$_{}$}
\label{5.10}
\begin{em}

Let $H \subseteq G$ be an open subgroup and $(\sigma ,  W)$ smooth. Then there is a functorial isomorphism
\[  {\rm Hom}_{G}(  c-{\rm Ind}_{H}^{G}(W) , \pi) \stackrel{\cong}{\longrightarrow}  
{\rm Hom}_{H}(  W  , \pi )  \]
given by $F \mapsto F \cdot  f $, the composition with the $H$-map $f$ of Lemma \ref{5.7A}.
\end{em}
\end{proposition}

\begin{example}{$c-\underline{{\rm Ind}}_{H}^{G}(\phi)$}
\label{5.11}
\begin{em}

Suppose that $\phi : H \longrightarrow k^{*}$ is a continuous character (i.e. a one-dimensional smooth representation of $H$).

Suppose that we are in a situation analogous to that of Example \ref{5.8}. Namely suppose that $H$ is open and closed, contains $Z(G)$,the centre of $G$, and is compact open modulo $Z(G)$. A basis for $k$ is given by $1 \in k^{*}$ and we have the function $f_{1} \in X_{c}$ given by
$f_{1}(h) = \phi(h) $ if $h \in H$ and $f_{1}(g) =0$ if $g \not\in H$.

If, following Lemma \ref{5.7A}, ${\mathcal G}$ is a set of coset representatives for 
$H \backslash G  $ then a $k$-basis for $c-\underline{{\rm Ind}}_{H}^{G}(\phi)$ is given by 
 \[   \{  g \cdot f_{1}  \ | \  g \in {\mathcal G} \}  . \]
 
 For $g \in G$ we have
\[ \begin{array}{ll}
(g \cdot f_{1})(x) = f_{1}(xg) & =  \left\{  \begin{array}{cc}
0 & {\rm if } \ xg \not\in H, \\
\\
\phi(xg)  & {\rm if } \ xg \in H, \\
\end{array} \right.   \\
\\
& =  \left\{  \begin{array}{cc}
0 & {\rm if } \  x  \not\in Hg^{-1}, \\
\\
\phi(xg)  & {\rm if } \ x \in Hg^{-1}. \\
\end{array} \right. 
\end{array} \]

Before going further let us introduce the presence of $(H, \phi)$ into the notation.
 \end{em}
\end{example}
\begin{definition}
\label{5.12}
\begin{em}

Let $H$ be a closed subgroup of $G$ containing the centre, $Z(G)$, which is compact open modulo $Z(G)$. Let $\phi : H \longrightarrow k^{*}$ be a continuous character of $H$. Denote by $X_{c}(H, \phi)$ the $k$-vector space of functions $f: G \longrightarrow k  $ such that

(i) \  $f(hg) = \phi(h)f(g)$ for all $h \in H, g \in G$,

(ii)  \  there is a compact open subgroup $K_{f} \subseteq G$ such that $f(gk) = f(g)$ for all $g \in G, k \in K_{f}$,

(ii)  $f$ is compactly supported modulo $H$.

As in \S\ref{5.6A}, the left action of $G$ on  $X_{c}(H, \phi)$ is given by $(g \cdot f)(x)= f(xg)$ and therefore
\[   \Sigma :  G  \longrightarrow  {\rm Aut}_{k}( X_{c}(H, \phi))  \]
gives a smooth representation of $G$ - denoted by $\Sigma = c-{\rm Ind}_{H}^{G}(\phi)$.

Henceforth we shall denote the map written as $f_{1}$ in Example \ref{5.11} by 
\linebreak
$f_{(H, \phi)} \in 
X_{c}(H, \phi)$.

Therefore, for $g \in G$, we have
\[ \begin{array}{ll}
(g \cdot f_{(H, \phi)})(x) = f_{(H, \phi)}(xg) & =  \left\{  \begin{array}{cc}
0 & {\rm if } \ xg \not\in H, \\
\\
\phi(xg)  & {\rm if } \ xg \in H, \\
\end{array} \right.   \\
\\
& =  \left\{  \begin{array}{cc}
0 & {\rm if } \  x  \not\in Hg^{-1}, \\
\\
\phi(xg)  & {\rm if } \ x \in Hg^{-1}. \\
\end{array} \right. 
\end{array} \]
\end{em}
\end{definition}

\begin{dummy}{{\rm Trace of} $c-{\rm Ind}_{H}^{GL_{n}F}(k_{\phi}) $}
\label{A5.13}
\begin{em}

Let $K_{0} = {\rm Ker}(\phi : H \longrightarrow k^{*} )$. We have $F^{*} \subset H$ and the valuation extension
\[ 0 \longrightarrow {\mathcal O}_{F}^{*} \longrightarrow F^{*} \longrightarrow {\mathbb Z} \longrightarrow 0 \]
in which ${\mathcal O}_{F}^{*} $ is compact open. Since $H$ is compact open modulo $F^{*}$ the subgroup $K_{0} \subset GL_{n}F$ is compact, open.

Consider $c-{\rm Ind}_{H}^{GL_{n}F}(k_{\phi})^{K_{0}} $. If $g \in K_{0}$ then 
\[ \begin{array}{l}
(g \cdot f_{(H, \phi)})(x) = f_{(H, \phi)}(xg) \\
\\
 =  \left\{  \begin{array}{cc}
0 & {\rm if } \  x  \not\in Hg^{-1}, \\
\\
\phi(xg)  & {\rm if } \ x \in Hg^{-1}. \\
\end{array} \right. \\
\\
 =  \left\{  \begin{array}{cc}
0 & {\rm if } \  x  \not\in Hg^{-1}, \\
\\
\phi(x)  & {\rm if } \ x \in Hg^{-1}. \\
\end{array} \right. \\
\\
=  f_{(H, \phi)}(x)
\end{array} \]
so that $ f_{(H, \phi)} \in c-{\rm Ind}_{H}^{GL_{n}F}(k_{\phi})^{K_{0}} $.
\end{em}
\end{dummy}

We have a chain of inclusions $K_{0} \leq H \leq N_{GL_{n}F}(H, \phi) \leq N_{GL_{n}F}(K_{0})$ where 
$N_{GL_{n}F}$ denotes a normaliser subgroup. If $ g \in N_{GL_{n}F}(K_{0})$ and $g_{1} \in K_{0}$ then
\[ \begin{array}{l}
g_{1}((g \cdot f_{(H, \phi)}))(x)   \\
\\
= ( g_{1}g) \cdot f_{(H, \phi)})(x)   \\
\\
= f_{(H, \phi)}(xg_{1}g) \\
\\
 =  \left\{  \begin{array}{cc}
0 & {\rm if } \  xg_{1}g  \not\in H, \\
\\
\phi(xg_{1}g)  & {\rm if } \ xg_{1}g \in H. \\
\end{array} \right. \\
\\
 =  \left\{  \begin{array}{cc}
0 & {\rm if } \  xg_{1}g  \not\in H, \\
\\
\phi(xg)  \phi(g^{-1}g_{1}g)  & {\rm if } \ xg_{1}g \in H. \\
\end{array} \right. \\
\\
 =  \left\{  \begin{array}{cc}
0 & {\rm if } \  xg  \not\in H, \\
\\
\phi(xg)    & {\rm if } \ xg \in H. \\
\end{array} \right. \\
\\
=  (g \cdot f_{(H, \phi)})(x)
\end{array} \]
so that, if $g \in N_{GL_{n}F}(K_{0})$ then $g \cdot f_{(H, \phi)} \in  c-{\rm Ind}_{H}^{GL_{n}F}(k_{\phi})^{K_{0}} $.

Since ${\rm dim}_{k}(c-{\rm Ind}_{H}^{GL_{n}F}(k_{\phi})^{K_{0}} ) < \infty$ we must have 
$[  N_{GL_{n}F}(K_{0}) : H ] < \infty$.

\section{The trace of  $c-{\rm Ind}_{H}^{GL_{n}F}(k_{\phi})$}

From (\cite{Sn20} \S 6 and \S13) we recapitulate the definition of the trace of an admissible $(\pi, V)$ of $G = GL_{n}F$ such as the case when $(\pi, V)$ is given by the representation on $V =  c-{\rm Ind}_{H}^{GL_{n}F}(k_{\phi})$.

Let ${\mathcal H}_{G}$ be, as before,  the space of smooth compactly supported complex-valued functions on $X=G$. Assuming $G$ is unimodular ${\mathcal H}_{G}$ is an algebra without unit under the convolution product
\[    (\phi_{1} * \phi_{2})(g) = \int_{G} \ \phi_{1}(gh^{-1}) \phi_{2}(h) dh .  \]
This is the Hecke algebra - an idempotented algebra (see \cite{Sn20} \S6). This is the case when $G = GL_{n}F$.

If $\phi \in {\mathcal H}_{GL_{n}F}$ define $\pi(\phi) \in {\rm End}(V)$ with $V$ as above 
\[   \pi(\phi)(v) = \int_{G} \phi(g) \pi(g)(v) dg  . \]
Then 
\[    \pi (\phi_{1} * \phi_{2})  =  \pi(\phi_{1}) \cdot \pi(\phi_{2})  \]
so that $V$ is an ${\mathcal H}_{GL_{n}F}$-representation.

The integral defining $\phi$ may be replaced by a finite sum as follows. Choose an open subgroup $K_{0}$ fixing $v$. Choosing $K_{0}$ small enough we may assume that the support of $\phi$ is contained in a finite union of left cosets $\{ g_{i}K_{0} \ | \ 1 \leq i \leq t \}$. Then
\[   \pi(\phi)(v) =  \frac{1}{{\rm vol}(K_{0})}  \sum_{i=1}^{t}  \phi(g_{i})  \pi(g_{i})(v) . \]

As with representations of finite groups the character of an admissible representation of a totally disconnected locally compact group $G$ is an important invariant. It is a distribution. It is a theorem of Harish-Chandra that if $G$ is a reductive $p$-adic group then the character is in fact a locally integrable function defined on a dense subset of $G$.

We shall define the character as a distribution on ${\mathcal H}_{G} = C_{c}^{\infty}(G)$. Suppose that $U$ is a finite-dimensional vector space and let $f:U \longrightarrow U$ be a linear map. Suppose ${\rm Im}(f) \subseteq U_{0} \subseteq U$. Then we have 
\[ {\rm Trace}(f:U_{0} \longrightarrow U_{0}) = {\rm Trace}(f:U \longrightarrow U) . \]
Therefore we may define the trace of any endomorphism $f$ of $V$ which has finite rank by choosing any
finite-dimensional $U_{0}$ such that ${\rm Im}(f) \subseteq U_{0} \subseteq V$ and by defining
\[   {\rm Trace}(f) =   {\rm Trace}(f:U_{0} \longrightarrow U_{0}) . \]

Now let $(\pi, V)$ be an admissible representation of $G$. Let $\phi \in {\mathcal H}_{G}$. Since $\phi$ is compactly supported and locally constant there exists a compact open $K_{0}$ such that $\phi \in {\mathcal H}_{K_{0}} = e_{K_{0}} *  {\mathcal H}_{G} * e_{K_{0}} $. Here $ e_{K_{0}} \in {\mathcal H}_{G} $ is the idempotent given by
\[ e_{K_{0}} = \frac{1}{{\rm vol}(K_{0})} \chi_{K_{0}}  \]
where $ \chi_{K_{0}} $ is the characteristic function of $K_{0}$.

The endomorphism $\pi(\phi)$ has image in $V^{K_{0}}$ which is finite-dimensional - by admissibility - so we define the trace distribution
\[ \chi_{V} : {\mathcal H}_{G} \longrightarrow {\mathbb C} \]
by
\[   \chi_{V}(\phi) = {\rm Trace}(\pi(\phi)) ,\]
which is called the trace of $V$.

Let us pause to consider the condition $\phi \in  e_{K_{0}} *  {\mathcal H}_{G} * e_{K_{0}} $.
\[ \begin{array}{l}
 e_{K_{0}} *  \phi * e_{K_{0}}(g) \\
 \\
 = \int_{G}   (e_{K_{0}} *  \phi )(gh^{-1}) \frac{1}{{\rm vol}(K_{0})} \chi_{K_{0}}(h) dh \\
 \\
  = \int_{K_{0}}   (e_{K_{0}} *  \phi )(gh^{-1}) \frac{1}{{\rm vol}(K_{0})}  dh \\
 \\
 =  \int_{h \in K_{0}}  \int_{G} \frac{1}{{\rm vol}(K_{0})} \chi_{K_{0}}( gh^{-1} h_{1}^{-1}) \phi(h_{1}) dh_{1}  \frac{1}{{\rm vol}(K_{0})}  dh \\
 \\
  =  \int_{h \in K_{0}}  \int_{gh^{-1}h_{1}^{-1} \in K_{0} } \chi_{K_{0}}( gh^{-1} h_{1}^{-1}) \phi(h_{1}) dh_{1}  \frac{1}{{\rm vol}(K_{0})^{2}}  dh \\
 \\
\end{array} \]
Since $h_{1}^{-1} \in K_{0}gK_{0}$ in the integral we may, by shrinking $K_{0}$ if necessary, ensure that both $\phi$ and $\phi((-)^{-1})$ are constant on $ K_{0}gK_{0}$ so that $\phi(h_{1}) = \phi(g)$. In fact we can ensure that this is true for all $ K_{0}gK_{0}$ because $\phi$ is compactly supported and locally constant.
Hence for all $g \in G$ we have  $e_{K_{0}} *  \phi * e_{K_{0}}(g) = \phi(g)$ and so 
$\phi \in e_{K_{0}} *  {\mathcal H}_{G} *  e_{K_{0}} = {\mathcal H}_{K_{0}} $.

In order to calculate the trace of $V = c-{\rm Ind}_{H}^{GL_{n}F}(k_{\phi})$. By Lemma \ref{5.7A} (iii) a $k$-basis for $V$ is $\{ g_{s} \cdot f_{(H, \phi)}  \ |  \ s \geq 1 \}$ where the $g_{i}$'s are a set of coset representatives for $G/H$.

If $\lambda \in {\mathcal H}_{GL_{n}F}$ we shall suppose that $K'_{0}$ is compact open and such that the support of $\lambda$ is contained in a finite union of left cosets $z_{i}K'_{0}$ so that we have a formula of the form
\[  \begin{array}{l}
 \pi(\lambda)(g_{s} \cdot f_{(H, \phi)}) \\
 \\
 =  \frac{1}{{\rm vol}(K'_{0})}  \sum_{i=1}^{t}  \lambda(z_{i})   z_{i} g_{s} \cdot f_{(H, \phi)} . 
 \end{array}  \]
for $ s \geq 1$. 

It would seem to be a good idea to choose $K'_{0} \leq K_{0} = {\rm Ker}(\phi : H \longrightarrow k^{*})$. Then we have  $K'_{0} \leq K_{0} \leq H \leq N_{GL_{n}F}K_{0}$ and
\[  \begin{array}{l}
  \chi_{ c-{\rm Ind}_{H}^{GL_{n}F}(k_{\phi})}(\lambda) \\
  \\
  = {\rm Trace}(\pi(\lambda))  \\
\\
 =  \frac{1}{{\rm vol}(K'_{0})}  \sum_{i=1}^{t}  \sum_{s \geq 1,   z_{i} g_{s} \cdot f_{(H, \phi)} = \mu_{i,s}  f_{(H, \phi)} } \lambda(z_{i})   \mu_{i,s} .
  \end{array}  \]

  Suppose that $ z_{i} g_{s} \cdot f_{(H, \phi)} = \mu_{i,s}  f_{(H, \phi)} $ then
  \[ \begin{array}{l}
   z_{i} g_{s} \cdot f_{(H, \phi)}(x) \\
   \\
   =  f_{(H, \phi)}(x  z_{i} g_{s} ) \\
\\
 =  \left\{  \begin{array}{cc}
0 & {\rm if } \  x  z_{i} g_{s}  \not\in H, \\
\\
\phi(x  z_{i} g_{s})  & {\rm if } \ x  z_{i} g_{s} \in H. \\
\end{array} \right. \\
\\
= \mu_{i,s}  f_{(H, \phi)}(x ) \\
\\
 =  \left\{  \begin{array}{cc}
0 & {\rm if } \  x   \not\in H, \\
\\
\mu_{i,s}  \phi(x)   & {\rm if } \ x   \in H. \\
\end{array} \right. \\
\\
  \end{array} \]
  which implies that $ z_{i} g_{s} \in H$ and $\mu_{i,s} = \phi(z_{i}g_{s})$. Since $g_{s}$ is one of the coset representatives in $G/H$ we must have $z_{i}g_{s} = h_{i,s} \in H$ and $g_{s} h_{i,s} ^{-1}  = z_{1}^{-1} $.
  Therefore $\phi(h_{i,s}) z_{i}^{-1} f_{(H,\phi)} = g_{s} \cdot f_{(H, \phi)}$, which is easily verified.
  
    For a fixed $\lambda$ there are only a finite number $t$ of $z_{1}$'s each defining a coset $z_{i}K'_{0}$ which may be assumed distinct. Therefore the formula for the trace 
  \[  \begin{array}{l}
  \chi_{ c-{\rm Ind}_{H}^{GL_{n}F}(k_{\phi})}(\lambda) \\
  \\
  = {\rm Trace}(\pi(\lambda))
\\
\\
 =  \frac{1}{{\rm vol}(K'_{0})}  \sum_{i=1}^{t}  \sum_{s \geq 1,   z_{i} g_{s} \in H}  \lambda(z_{i})   \phi(z_{i}g_{s}) 
  \end{array}  \]
  makes sense as a finite sum since given $z_{i}$ only a finite number of $g_{s}$'s satisfy the condition in the sum, as $K'_{0} \leq H$.

{\it DS: At this point my father inserted the comment `GOT TO HERE'. I don't know whether he simply omitted removing this comment or if there is in fact something missing here.}

\begin{definition}
\label{5.13}
\begin{em}

For $(H, \phi)$ and $(K, \psi)$ as in Definition \ref{5.12},  write 
\linebreak
$[(K, \psi), g, (H, \phi)]$ for any triple consisting of $g \in G$, characters $\phi, \psi$ on 
subgroups $H, K \leq G$, respectively such that 
\[    (K, \psi) \leq (g^{-1}Hg, (g)^{*}(\phi)) \]
which means that $K \leq  g^{-1}Hg$ and that $\psi(k) = \phi(h)$ where $k = g^{-1}hg$ for $h \in H, k \in K$.

Let ${\mathcal H}$ denote the $k$-vector space with basis given by these triples. Define a product on these triples by the formula 
\[  [(H, \phi), g_{1}, (J, \mu)]  \cdot  [(K, \psi), g_{2}, (H, \phi)] =   [(K, \psi), g_{1}g_{2}, (J, \mu)]  \]
and zero otherwise. This product makes sense because 

(i)  \   if $K \leq g_{2}^{-1} H g_{2}$ and
 $H \leq g_{1}^{-1} J g_{1}$ then $K \leq  g_{2}^{-1} H g_{2} \leq   g_{2}^{-1} g_{1}^{-1} J g_{1} g_{2} $ 
 
 and 
 
 (ii)  \  if $\psi(k) = \phi(h) = \mu(j), $ where $k = g_{2}^{-1}hg_{2},  h = g_{1}^{-1}j g_{1}$ then
 \linebreak
  $k = g_{2}^{-1}  g_{1}^{-1}j g_{1}   g_{2}$. 
 
 This product is clearly associative and we define an algebra ${\mathcal H}_{cmc}(G)$ to be ${\mathcal H}$ modulo the relations (c.f. Appendix \S4)
 \[   [(K, \psi), gk, (H, \phi)]  = \psi(k^{-1}) [(K, \psi), g, (H, \phi)]  \]
and
 \[     [(K, \psi), hg, (H, \phi)]  = \phi(h^{-1}) [(K, \psi), g, (H, \phi)] .    \]
 
 We observe that
 \[  [(K, \psi), g, (H, \phi)] = [( g^{-1}Hg, g^{*}\phi), g , (H, \phi)] \cdot [(K, \psi), 1, ( g^{-1}Hg, g^{*}\phi)]\]
 
We shall refer to this algebra as the compactly supported modulo the centre (CSMC-algebra) of $G$.
\end{em}
\end{definition}

\begin{lemma}{$_{}$}
\label{5.14}
\begin{em}

Let $[(K, \psi), g, (H, \phi)]$ be a triple as in Definition \ref{5.13}. Associated to this triple define a left $k[G]$-homomorphism
\[  [(K, \psi), g, (H, \phi)] : X_{c}(K, \psi) \longrightarrow X_{c}(H, \phi) \]
by the formula $g_{1} \cdot f_{(K, \psi)} \mapsto (g_{1}g^{-1}) \cdot  f_{(H, \phi) }$.
\end{em}
\end{lemma}

For a proof, which is the same as in the case when $G$ is finite, can be found in (the Appendix on induction in the case of finite groups).

\begin{theorem}{$_{}$}
\begin{em}
\label{5.15}

Let ${\mathcal M}_{c}(G)$ denote the partially order set of pairs $(H, \phi)$ as in Definitions \ref{5.12} and \ref{5.13} (so that $ X_{c}(H, \phi) = c-{\rm Ind}_{H}^{G}(\phi)$). Then, when each $n_{\alpha} =1$, 
\[        M_{c}(\underline{n}, G) = \oplus_{\alpha \in {\mathcal A}, (H, \phi) \in {\mathcal M}_{c}(G)}  \ n_{\alpha} X_{c}(H, \phi)  \]
is a left $k[G] \times {\mathcal H}_{cmc}(G)$-module. For a general distribution of multiplicities
$ \{ n_{\alpha} \}$ it is Morita equivalent to a left $k[G] \times {\mathcal H}_{cmc}(G)$-module.
\end{em}
\end{theorem}

{\bf Proof}

We have only to verify associativity of the module multiplication, which is obvious.  $\Box$

\begin{definition}{$_{k[G]}{\bf mon}$, the monomial category of $G$}
\label{5.16}
\begin{em}

The monomial category of $G$ is the additive category (non-abelian) whose objects are the $k$-vector spaces given by direct sums of $ X_{c}(H, \phi)$'s of \S\ref{5.15} and whose morphisms are elements of the hyperHecke algebra ${\mathcal H}_{cmc}(G)$. In other words the subcategory of the category of 
$k[G] \times {\mathcal H}_{cmc}(G)$-modules of which one example is $ M_{c}(\underline{n}, G)$ in \S\ref{5.15}.
\end{em}
\end{definition}

{\it DS: At this point my father inserted the comment `THIS IS SECTION 13 FROM Derived Langlands II' and the Appendix below is the Appendix from Section 13 of Derived Langlands II. }

\section{Appendix: Smooth representations and Hecke modules}

In this Appendix, for my convenience, representations are complex representations.

Now let $\Gamma$ be a compact totally disconnected group. Denote by $\hat{\Gamma}$ the set of equivalence classes of finite-dimensional irreducible representations of $\Gamma$ whose kernel is open - and hence of finite index in $\Gamma$.

Suppose now that $\Gamma$ is finite  and $(\pi, V)$ is a representation of $\Gamma$ on a possible infinite dimensional vector space $V$. If $\rho \in \hat{\Gamma}$ let $V(\rho)$ be the sum of all invariant subspaces of $V$ that are isomorphic as $\Gamma$-modules to $V_{\rho}$. $V(\rho)$ is the $\rho$-isotypic subspace of $V$. We have
\[ V \cong \oplus_{\rho \in \hat{\Gamma}} \ V_{\rho} . \]

Now we generalise this to smooth representations of a totally disconnected locally compact group. Choose a compact open subgroup $K$ of $G$. The compact open normal subgroups of $K$ form a basis of neighbourhoods of the identity in $K$. Let $\rho  \in \hat{K}$ then the kernel of $\rho$ is $K_{\rho}$ a compact open normal subgroup of finite index.
\begin{proposition}{(\cite{DB96} Proposition 4.2.2)}
\label{12.2}
\begin{em}

Let $(\pi, V)$ be a smooth representation of $G$.
Then
\[ V \cong \oplus_{\rho \in \hat{K}} \ V_{\rho} . \]
The representation $\pi$ is admissible if and only if each $V(\rho)$ is finite-dimensional.
\end{em}
\end{proposition}

Let $(\pi, V)$ be a smooth representation of $G$. If $\hat{v} : V \longrightarrow  {\mathbb C}$ is a linear functional we write $\langle v , \hat{v} \rangle = \hat{v}(v)$ for $v \in V$. We say $\hat{v}$ is smooth if there exists an open neighbourhood $U$ of $1 \in G$ such that for all $g \in U$
\[  \langle \pi(g)(v) , \hat{v} \rangle = \hat{v}(v) . \]
Let $\hat{V}$ denote the space of smooth linear functionals on $V$.

Define the contragredient representation $(\hat{\pi}, \hat{V})$ is defined by
\[   \langle v ,   \hat{\pi}(g)(\hat{v}) \rangle =   \langle \pi(g^{-1})(v) , \hat{v} \rangle . \]
The contragredient representation of a smooth representation is a smooth representation. Also
\[ \hat{V} \cong \oplus_{\rho \in \hat{K}} \ V_{\rho}^{*}  \]
where $V_{\rho}^{*}$ is the dual space of $V_{\rho}$.

Since the dual of a finite-dimensional $V_{\rho}$ is again finite-dimensional the contragredient of an admissible representation is also admissible. Also $\hat{\hat{\pi}} = \pi$.

If $X$ is a totally disconnected space a complex valued function $f$ on $X$ is smooth if it is locally constant.
Let ${\mathcal H}_{G}$ be, as before,  the space of smooth compactly supported complex-valued functions on $X=G$. Assuming $G$ is unimodular ${\mathcal H}_{G}$ is an algebra without unit under the convolution product
\[    (\phi_{1} * \phi_{2})(g) = \int_{G} \ \phi_{1}(gh^{-1}) \phi_{2})(h) dh .  \]
This is the Hecke algebra - an idempotented algebra (see \S6).

If $\phi \in {\mathcal H}$ define $\pi(\phi) \in {\rm End}(V)$ with $V$ as above 
\[   \pi(\phi)(v) = \int_{G} \phi(g) \pi(g)(v) dg  . \]
Then 
\[    \pi (\phi_{1} * \phi_{2})  =  \pi(\phi_{1}) \cdot \pi(\phi_{2})  \]
so that $V$ is an ${\mathcal H}$-representation.

The integral defining $\phi$ may be replaced by a finite sum as follows. Choose an open subgroup $K_{0}$ fixing $v$. Choosing $K_{0}$ small enough we may assume that the support of $\phi$ is contained in a finite union of left cosets $\{ g_{i}K_{0} \ | \ 1 \leq i \leq t \}$. Then
\[   \pi(\phi)(v) =  \frac{1}{{\rm vol}(K_{0})}  \sum_{i=1}^{t}  \phi(g_{i})  \pi(g_{i})(v) . \]

{\bf Finite group example:}

Let $(\pi, V)$ be a finite-dimensional representation of a finite group $G$. Write ${\mathcal H}$ for the space of functions from $G$ to ${\mathbb C}$. If $\phi_{1}, \phi_{2} \in {\mathcal H}$ define $\phi_{1} * \phi_{2} 
\in {\mathcal H}$ by
\[   (\phi_{1} * \phi_{2})(g) = \sum_{h \in G} \ \phi_{1}(gh^{-1}) \phi_{2}(h) . \]
For $\phi \in {\mathcal H}$ define $\pi(\phi) \in {\rm End}_{{\mathbb C}}(V)$ by
\[   \pi(\phi)(v) = \sum_{g \in G} \ \phi(g) \pi(g)(v) . \]
Hence
\[ \begin{array}{l}
\pi(\phi_{1}(  \pi(\phi_{2})(v) )   \\
\\
=  \pi(\phi_{1})(   \sum_{g \in G} \ \phi_{2}(g) \pi(g)(v)  ) \\
\\
=   \sum_{g \in G} \ \phi_{2}(g)  \pi(\phi_{1}( \pi(g)(v))  \\
\\
=    \sum_{g \in G} \   \phi_{2}(g)  \sum_{\tilde{g} \in G} \ \phi_{1}(\tilde{g}) (\pi(\tilde{g}( \pi(g)(v)) \\
\\
=   \sum_{g, \tilde{g} \in G} \  \phi_{2}(g)  \phi_{1}(\tilde{g}) (\pi(\tilde{g}g)(v)) .
\end{array} \]

Now 
\[ \begin{array}{l}
\pi( \phi_{1} * \phi_{2})(v) \\
\\
= \sum_{g_{1} \in G} \ ( \phi_{1} * \phi_{2})(g_{1}) \pi(g_{1})(v)  \\
\\
= \sum_{g_{1}, h \in G} \  \phi_{1}(h_{1}h^{-1}) \phi_{2}(h) \pi(g_{1})(v)  .
\end{array} \]
Setting $g=h$, $\tilde{g}g = g_{1}$ shows that 
\[  \pi( \phi_{1} * \phi_{2}) =  \pi(\phi_{1}) \cdot \pi(\phi_{2}) .\]

Also ${\mathcal H} \cong {\mathbb C}[G]$  because if $f_{g}(x) = 0$ if $g \not= x$ and $f_{g}(g)=1$ then 
\[ f_{g} * f_{g'} = f_{gg'} . \]

\begin{proposition}{(\cite{DB96} Proposition 4.2.3)}
\label{12.3}
\begin{em}

Let $(\pi, V)$ be a smooth non-zero representation of $G$. Then equivalent are:

(i)  \ $\pi$ is irreducible.

(ii) \ $V$ is a simple ${\mathcal H}$-module.

(iii) \ $V^{K_{0}}$ is either zero or simple as an ${\mathcal H}_{K_{0}}$-module for all open subgroups $K_{0}$.
Here ${\mathcal H}_{K_{0}} = e_{K_{0}} * {\mathcal H}  * e_{K_{0}}$.
\end{em}
\end{proposition}

Schur's Lemma holds (\cite{DB96} \S4.2.4)for $(\pi, V)$ an irreducible admissible representation of a totally disconnected group $G$.
\begin{proposition}{(\cite{DB96} Proposition 4.2.5)}
\label{12.4}
\begin{em}

Let $(\pi, V)$ be an admissible representation of the totally disconnected locally compact group $G$ with contragredient $(\hat{\pi}, \hat{V})$. Let $K_{0} \subseteq G$ be a compact open subgroup. Then the canonical pairing between $V$ and $\hat{V}$ induces a non-degenerate pairing betweem $V^{K_{0}}$ and 
$\hat{V}^{K_{0}}$.
\end{em}
\end{proposition}

{\bf The trace}

As with representations of finite groups the character of an admissible representation of a totally disconnected locally compact group $G$ is an important invariant. It is a distribution. It is a theorem of Harish-Chandra that if $G$ is a reductive $p$-adic group then the character is in fact a locally integrable function defined on a dense subset of $G$.

We shall define the character as a distribution on ${\mathcal H}_{G} = C_{c}^{\infty}(G)$. Suppose that $U$ is a finite-dimensional vector space and let $f:U \longrightarrow U$ be a linear map. Suppose ${\rm Im}(f) \subseteq U_{0} \subseteq U$. Then we have 
\[ {\rm Trace}(f:U_{0} \longrightarrow U_{0}) = {\rm Trace}(f:U \longrightarrow U) . \]
Therefore we may define the trace of any endomorphism $f$ of $V$ which has finite rank by choosing any
finite-dimensional $U_{0}$ such that ${\rm Im}(f) \subseteq U_{0} \subseteq V$ and by defining
\[   {\rm Trace}(f) =   {\rm Trace}(f:U_{0} \longrightarrow U_{0}) . \]

Now let $(\pi, V)$ be an admissible representation of $G$. Let $\phi \in {\mathcal H}_{G}$. Since $\phi$ is compactly supported and locally constant there exists a compact open $K_{0}$ such that $\phi \in {\mathcal H}_{K_{0}}$. The endomorphism $\pi(\phi)$ has image in $V^{K_{0}}$ which is finite-dimensional - by admissibility - so we define the trace distribution
\[ \chi_{V} : {\mathcal H} \longrightarrow {\mathbb C} \]
by
\[   \chi_{V}(\phi) = {\rm Trace}(\pi(\phi)) . \]
\begin{proposition}{(\cite{DB96} Proposition 4.2.6)}
\label{12.5}
\begin{em}

Let $R$ be an algebra over a field $k$. Let $E_{1}$ and $E_{2}$ be simple $R$-modules that are finite-dimensional over $k$. For each $\phi \in R$ if 
\[ {\rm Trace}( (\phi \cdot -):E_{1} \longrightarrow E_{1}) =  {\rm Trace}( (\phi \cdot -):E_{2} \longrightarrow E_{2})   \]
then the $E_{i}$ are isomorphic $R$-modules.
\end{em}
\end{proposition} 
\begin{proposition}{(\cite{DB96} Proposition 4.2.7)}
\label{12.6}
\begin{em}

Let $(\pi_{1}, V_{1})$ and $(\pi_{2}, V_{2})$ be irreducible admissible representations of $G$ (as above) such that, for each compact open $K_{1}$, $V_{1}^{K_{1}} \cong V_{2}^{K_{1}}$ as ${\mathcal H}_{K_{1}}$-modules then $(\pi_{1}, V_{1}) \cong (\pi_{2}, V_{2})$.
\end{em}
\end{proposition} 
\begin{theorem}{(\cite{DB96} Theorem 4.2.1)}
\label{12.7}
\begin{em}

Let $(\pi_{1}, V_{1})$ and $(\pi_{2}, V_{2})$ be irreducible admissible representations of $G$ (as above) such that $\chi_{V_{1}} = \chi_{V_{2}}$ then $(\pi_{1}, V_{1}) \cong (\pi_{2}, V_{2})$.
\end{em}
\end{theorem} 

From this one sees that the contragredient of an admissible irreducible $(\pi, V)$ of $GL_{n}K$ ($K$ a $p$-adic local field) 
is given by $\pi_{1}(g) = \pi( (g^{-1})^{tr}) $
on the same vector space $V$.

 \end{document}